\documentclass[twoside,11pt]{article}
\title{} \author{} \date{}
\usepackage{latexsym,amssymb,times}
\input amssym.def
\newtheorem{te}{Theorem}[section]

\newtheorem{fac}[te]{Fact}
\newtheorem{cla}[te]{Claim}

\def\dok{\noindent{\bf Proof. }}
\def\kdok{\hfill $\Box$ \par \vspace*{2mm} }
\def\a{\alpha}
\def\f{\varphi}
\def\p{\psi}
\def\o{\omega}

\def\r{\rho}

\def\N{{\mathbb N}}
\def\X{{\mathbb X}}
\def\Y{{\mathbb Y}}
\def\Z{{\mathbb Z}}
\def\A{{\mathbb A}}
\def\BL{{\mathbb L}}
\def\BH{{\mathbb H}}
\def\BK{{\mathbb K}}
\def\BI{{\mathbb I}}
\def\BF{{\mathbb F}}
\def\BM{{\mathbb M}}
\def\CT{{\mathcal T}}
\def\CK{{\mathcal K}}
\def\L{{\mathcal L}}
\def\CM{{\mathcal M}}
\def\c{{\mathfrak{c}}}
\def\la{\langle}
\def\ra{\rangle}
\def\id{\mathop{\mathrm{id}}\nolimits}
\def\Iso{\mathop{\rm Iso}\nolimits}
\def\Aut{\mathop{\rm Aut}\nolimits}
\def\Sym{\mathop{\rm Sym}\nolimits}
\def\Sent{\mathop{\rm Sent}\nolimits}
\def\Form{\mathop{\rm Form}\nolimits}
\def\Age{\mathop{\mathrm{Age}}\nolimits}
\def\ar{\mathop{\rm ar}\nolimits}
\def\Mod{\mathop{\rm Mod}\nolimits}
\def\LO{\mathop{\rm LO}\nolimits}
\def\Th{\mathop{\rm Th}\nolimits}
\def\Lit{\mathop{\rm Lit}\nolimits}
\def\Pa{\mathop{\rm Pa}\nolimits}
\def\dom{\mathop{\rm dom}\nolimits}
\def\ran{\mathop{\rm ran}\nolimits}
\def\Ker{\mathop{\rm Ker}\nolimits}
\begin{document}
\thispagestyle{plain}
\begin{center}
           {\large \bf \uppercase{Vaught's Conjecture for Almost Chainable  Theories}}
\end{center}
\begin{center}
{\bf Milo\v s S.\ Kurili\'c}\footnote{Department of Mathematics and Informatics, Faculty of Sciences, University of Novi Sad,
                                      Trg Dositeja Obradovi\'ca 4, 21000 Novi Sad, Serbia.
                                      e-mail: milos@dmi.uns.ac.rs}
\end{center}
\begin{abstract}
\noindent
A  structure $\Y$ of a relational language $L$ is called {\it almost chainable}
iff there are a finite set $F \subset Y$ and a linear order $\,<$ on the set $Y\setminus F$
such that for each partial automorphism $\f$ (i.e., local automorphism, in Fra\"{\i}ss\'{e}'s terminology)
of the linear order $\la Y\setminus F, <\ra$
the mapping $\id _F \cup \f$ is a partial automorphism of $\Y$.
By a theorem of Fra\"{\i}ss\'{e}, if $|L|<\o$, then $\Y$ is almost chainable iff the profile of $\Y$ is bounded;
namely, iff there is a positive integer $m$ such that $\Y$ has $\leq m$ non-isomorphic substructures of size $n$, for each positive integer $n$.
A complete first order $L$-theory $\CT$ having infinite models is called almost chainable
iff all models of $\CT$ are almost chainable
and it is shown that the last condition is equivalent to the existence of one countable almost chainable model of $\CT$.
In addition, it is proved  that
an almost chainable theory has either one or continuum many non-isomorphic countable models
and, thus,  the Vaught conjecture is confirmed for almost chainable  theories.


\end{abstract}
\section{Introduction}\label{S1}
In this article we confirm Vaught's conjecture
for almost chainable theories,
extending the result of \cite{K}, which concerns the smaller class of monomorphic theories.
We recall that the Vaught conjecture is related to the number $I(\CT ,\o)$ of non-isomorphic countable models
of a countable complete first order theory $\CT$.
In 1959 Robert Vaught \cite{Vau} asked if there is a theory $\CT$
such that the equality $I(\CT ,\o)=\o _1$ is provable without the use of the continuum hypothesis;
since then, the implication $I(\CT ,\o)> \o \Rightarrow I(\CT ,\o)= \c$ is known as Vaught's conjecture.

The rich history of the investigation related to that (still unresolved) conjecture
includes a long list of results confirming the conjecture in particular classes of theories
(see, for example, the introduction and references of \cite{Pun}) and, on the other hand,
intriguing results concerning the consequences of the existence of counterexamples and the properties of (potential) counterexamples
(see, e.g., \cite{Bal}).

The results of this paper are built on the fundament consisting of two (groups of) results.
The first one is the basic Rubin's paper \cite{Rub} from 1974 in which the Vaught conjecture is confirmed
for theories of linear orders with unary predicates;
we will use the following result of Rubin  (see Theorem 6.12 of \cite{Rub}).
\begin{te}[Rubin]\label{T002}
If $\CT$ is a complete theory of a linear order with a finite set of unary predicates, then $I(\CT, \o )\in \{ 1, \c\}$.
\end{te}
The second group of results is a part of the fundamental work
concerning combinatorial properties of relational structures collected in the book of Roland Fra\"{\i}ss\'{e} \cite{Fra}.
We will use  Fra\"{\i}ss\'{e}'s  results related to almost chainable structures,
as well as a theorem of Gibson, Pouzet and Woodrow from \cite{Gib},
describing all linear orders which chain an almost chainable structure over a fixed finite set,
which is derived from similar results obtained independently by Frasnay \cite{Fras}
and by Hodges, Lachlan and Shelah \cite{Hodg1}. These results are presented in Section \ref{S2}.

In Section \ref{S3} we show that
a complete theory $\CT $ with infinite models has a countable almost chainable model iff  all models of $\CT$ are almost chainable and, so,
establish the notion of an almost chainable theory.  In Section
\ref{S4} we prove that for each complete almost chainable theory $\CT$ having infinite models we have $I(\CT ,\o)\in \{ 1 , \c \}$ and,
thus, confirm the Vaught conjecture for such theories.

The results of this paper generalize the results of \cite{K} about theories of monomorphic structures\footnote{A relational structure $\Y$ is
monomorphic iff all its $n$-element substructures are isomorphic, for each positive integer $n$,
while (for $|L|<\o$, see \cite{Fra}, p.\ 297) $\Y$ is almost chainable iff there is a positive integer $m$ such that $\Y$ has $\leq m$ non-isomorphic substructures of size $n$, for each $n\in \N$.
}
and we note that the arguments used in our proofs are, as in \cite{K}, more combinatorial than model-theoretical.
Also we remark that some parts of this paper are (more or less) folklore or similar to the corresponding parts of \cite{K},
but, for completeness of the paper, they are included in the text.

\section{Preliminaries. Almost chainable structures}\label{S2}
Throughout the paper
we assume that $L=\la R_i :i\in I\ra$ is a relational language, where  $\ar (R_i)=n_i\in \N$, for $i\in I$.
If $Y$ is a non-empty set and $\CT \subset \Sent _L$ an $L$-theory, then
$\Mod _{L}^{\CT }(Y)$ (resp.\ $\Mod _{L}(Y)$; $\Mod _{L}^{\CT }$) will denote
the set of all models of $\CT$ with domain $Y$
(resp.\ the set of all $L$-structures with domain $Y$;
the class of all models of $\CT$).
Let $\Y =\la Y, \la R_i ^\Y : i\in I\ra\ra$ be an $L$-structure.
For a non-empty set $H\subset Y$, $\BH :=\la H, \la R_i ^\Y \upharpoonright  H : i\in I\ra\ra$ is the corresponding {\it substructure} of $\Y$.
If $J\subset I$, then $L_J:=\la R_i : i\in J\ra$  is the corresponding {\it reduction} of $L$
and $\Y |L_J :=\la Y, \la R_i ^\Y : i\in J\ra\ra$ the corresponding {\it reduct} of $\Y$.
By $[\Y]$ we will denote the class of all $L$-structures being isomorphic to $\Y$ (the isomorphism type of $\Y$).

If $\X=\la X ,<\ra$ is a linear order, then $\X ^*$ will denote its reverse, $\la X ,< ^{-1}\ra $.
By $LO _X$ we denote the set of all linear orders on the set $X$.

We recall the notions and concepts introduced by Fra\"{\i}ss\'{e} which will be used in this paper
and fix a convenient notation.
For $n\in \N$, by $\Age _n (\Y )$ we denote the collection $\{ [\BH ]:H\in [Y]^n \}$ of isomorphism types of $n$-element substructures of $\Y$
(or equivalently, $\Age _n (\Y )=\{ \BH \in \Mod _L (n): \BH\hookrightarrow \Y \}/\!\cong$).
The {\it age} of $\Y $ is the collection $\Age (\Y ):=\bigcup _{n\in \N }\Age _n (\Y )$.
The function $\f _\Y $ with the domain $\N$
defined by $\f _\Y (n)=|\Age _n (\Y )|$, for all $n\in \N$,  is the {\it profile} of $\Y$.

By $\Pa(\Y)$ we denote the set of all partial automorphisms of $\Y$ (isomorphisms between substructures of $\Y$,
or, in Fra\"{\i}ss\'{e}'s terminology, {\it local automorphisms}).
The $L$-structure $\Y$ is {\it freely interpretable} in an $L'$-structure $\X$ having the same domain
iff $\Pa (\X )\subset \Pa (\Y)$. We will say that $\Y$ is {\it simply definable} in $\X$ iff each relation $R_i^\Y$
is definable by a quantifier free $L'$-formula in the structure $\X$.
\paragraph{Almost chainable structures}
Let $\Y \in \Mod _L (Y)$, $F\in [Y]^{<\o}$ and $<\in LO _{Y\setminus F}$.
Following Fra\"{\i}ss\'{e} (see \cite{Fra}, p.\ 294), the structure $\Y$ is called {\it $(F,<)$-chainable} iff
\begin{equation} \label{EQ8085}
\forall \f \in \Pa (\la Y\setminus F ,<\ra) \;\; \id _F \cup \f \in \Pa (\Y ).
\end{equation}
The structure $\Y$ is called {\it $F$-chainable} if it is $(F,<)$-chainable for some linear order  $<$ on $Y\setminus F$.
$\Y$ is called {\it  almost chainable} if it is $F$-chainable for some $F\in [Y]^{<\o}$.

The following four statements are proved in \cite{Fra} for $|L|=1$ and have straightforward generalizations for
arbitrary relational language $L$. So, these results of Fra\"{\i}ss\'{e} are cited and used in the paper in such, more general, form.

Generally, if $Y$ is a set, $F\in [Y]^{n}$, $<\;\in LO _{Y\setminus F}$ and $F=\{ a_0,\dots ,a_{n-1}\}$ is an enumeration of the elements of the set $F$,
we introduce the auxiliary language $L_n:=\la R,U_0,\dots ,U_{n-1}\ra$,
consisting of new relational symbols, where $\ar (R)=2$ and $\ar (U_j)=1$, for $j<n$,
and define the linear order $\lhd$ on the set $Y$ and the $L_n$-structure (in fact, the linear order with $n$ unary predicates) $\X $ by \\[-7mm]
\begin{itemize}\itemsep -1.5mm \itemindent 2mm
\item[(L1)]  $\lhd \upharpoonright (Y\setminus F)=\;<$,
\item[(L2)] $\la Y ,\lhd\ra= \{a_0\}+\dots +\{a_{n-1}\}+(Y\setminus F) $,
\item[(L3)] $\X :=\la Y,\lhd, \{a_0\},\dots ,\{a_{n-1}\}\ra$.
\end{itemize}
\begin{fac}\label{T8081}
Let $\Y$ be an $L$-structure, $F=\{ a_0,\dots ,a_{n-1}\}\in [Y]^{n}$ and $<\;\in LO _{Y\setminus F}$. If $\lhd$ and $\X$ are defined by (L1)--(L3),
then the following conditions are equivalent:\\[-7mm]
\begin{itemize}\itemsep -1.5mm
\item[(a)] $\Y$ is $(F,<)$-chainable,
\item[(b)] $\Y$ is freely interpretable in $\X$,  (that is, $\Pa (\X )\subset \Pa (\Y)$),
\item[(c)] $\Y$ is simply definable in $\X$.
\end{itemize}
\end{fac}
\dok
(a) $\Rightarrow$ (b).
Let (\ref{EQ8085}) hold and $f\in \Pa (\X )$.
Then, since $f$ preserves $U_j$'s,
for each $y\in \dom f$ and each $j<n$ we have: $y=a_j$ iff $f(y)=a_j$.
So, if $a_j\in \dom f$, then $f(a_j)=a_j$
and, hence, $f \upharpoonright (F\cap \dom f)=\id _{F\cap \dom f}$.
In addition, if $y\in (\dom f )\setminus F$, then $f(y)\not\in F$
and, hence, $f[(\dom f )\setminus F]\subset Y\setminus F$.
So, by (L1), $\f :=f \upharpoonright ((\dom f )\setminus F)\in \Pa (\la Y\setminus F ,<\ra) $
and by (\ref{EQ8085}) we have $g:=\id _F \cup \;\f \in \Pa (\Y )$.
Finally, since $f=g\upharpoonright \dom f$, it follows that $f\in \Pa (\Y )$.

(b) $\Rightarrow$ (c).
Let $\Pa (\X )\subset \Pa (\Y )$ and $i\in I$.
For $\bar y\in Y^{n_i}$, let $\varepsilon _{\bar y}(v_0, \dots , v_{n_i-1})$
be the conjunction of all $L_n$-literals (i.e.\ atomic formulas and their negations)
in variables $v_0, \dots , v_{n_i-1}$
which are satisfied in the $L_n$-structure $\X$ by $\bar y$,
that is, $\varepsilon _{\bar y}(\bar v):=\bigwedge \{ \eta \in \Lit _{L_n}(\bar v): \X \models \eta [\bar y]\}$.
Since $|L_n|<\o$ the binary relation $\sim $ on the set $Y^{n_i}$
defined by $\bar x \sim \bar y$ iff $\varepsilon _{\bar x}=\varepsilon _{\bar y}$
is an equivalence relation with finitely many equivalence classes, say $m$.
If $\varepsilon _1 (\bar v), \dots ,\varepsilon _m (\bar v)$ is the list of the corresponding formulas
and, for $k\leq m$,
$D_{\varepsilon _k}:= \{ \bar y \in Y^{n_i}: \X \models \varepsilon _k [\bar y]\}$,
then $\{ D_{\varepsilon _k}: k\leq m \}$ is a partition of the set $Y^{n_i}$.
So, if we show that for each $k\leq m$ we have
\begin{equation}\label{EQ8066}
D_{\varepsilon _k}\cap R_i^\Y \neq \emptyset \Rightarrow D_{\varepsilon _k} \subset  R_i^\Y ,
\end{equation}
then $ R_i^\Y =\bigcup _{k\in J}D_{\varepsilon _k}$, where $J:=\{ k\leq m : D_{\varepsilon _k}\cap  R_i^\Y \neq \emptyset\}$,
and the relation $ R_i^\Y $ is definable in $\X $ by the quantifier-free $L_n$-formula $\f _i (\bar v):=\bigvee _{k\in J}\varepsilon _k (\bar v)$.

So, if $\bar x \in D_{\varepsilon _k}\cap  R_i^\Y $ and $\bar y \in D_{\varepsilon _k}$,
then $\X \models \varepsilon _k [\bar x ]$ and $\X \models \varepsilon _k [\bar y ]$
and, hence, $p:=\{ \la x_r, y_r\ra : r < n_i\}\in \Pa (\X )\subset \Pa (\Y )\subset \Pa (\la Y, \r _i \ra)$.
Thus, since $\bar x \in  R_i^\Y  \upharpoonright \{x_r : r< n_i \}$,
we have $\bar y =p\bar x \in  R_i^\Y  \upharpoonright \{y_r : r < n_i \}$
and, hence, $\bar y\in  R_i^\Y $.
So, (\ref{EQ8066}) is proved and $ R_i^\Y =\{ \bar y \in Y^{n_i}: \X \models \f _i [\bar y]\}$.

(c) $\Rightarrow$ (a).
For $i\in I$, let $\f _i (v_0, \dots , v_{n_i-1})\in \Form _{L_n}$ be a $\Sigma _0$-formula such that
\begin{equation}\label{EQ8067}
\forall \bar y \in Y^{n_i} \;\;\Big(\bar y\in R_i^\Y \Leftrightarrow \X  \models \f _i [\bar y ]\Big) .
\end{equation}
For a proof of (\ref{EQ8085}) we take $\f \in \Pa (\la Y\setminus F ,<\ra)$
and show that $f:=\id _F \cup\; \f \in \Pa (\Y )$.
By (L1) we have $\f \in \Pa (\la Y ,\lhd\ra)$; by (L2), $f \in \Pa (\la Y ,\lhd\ra)$ and, since $f(a_j)=a_j$, for all $j<n$, we obtain
$f\in \Pa (\X )$.
So, for $K:=\dom f$ and $H:=\ran f$,
denoting by $\BK$ and $\BH$ the corresponding substructures of $\X $, we have $f\in \Iso (\BK ,\BH)$.
Now, for  $i\in I$ and $\bar y \in K ^{n_i}$ we have
$\bar y\in R_i ^\Y$
iff (by (\ref{EQ8067})) $\X  \models \f _i [\bar y ]$
iff (since $\f_i$ is a $\Sigma _0$-formula) $\BK\models \f _i [\bar y ]$
iff (since $f\in \Iso (\BK ,\BH)$) $\BH\models \f _i [f \bar y ]$
iff (since $\f_i$ is a $\Sigma _0$-formula) $\X\models \f _i [f \bar y ]$
iff (by (\ref{EQ8067})) $f\bar y\in R_i ^\Y$.
So $f\in \Pa (\la Y, \r _i \ra)$, for all $i\in I$;
thus $f\in \Pa (\Y)$ and (\ref{EQ8085}) is true.
\kdok
\begin{fac}\label{T8088}
If $\Y$ is an infinite almost chainable $L$-structure,
then there is a minimal finite set $F\subset Y$ such that $\Y$ is $F$-chainable (the {\em kernel of} $\Y$, in notation $\Ker (\Y)$).
\end{fac}
\dok
For $|L|=1$, this is 10.9.3 of \cite{Fra}, p.\ 296. But the proof of 10.9.3 as well as the proofs of propositions (1), (2) and (3) of 10.9.2,
which are used in the proof of 10.9.3 have straightforward generalizations for arbitrary relational language $L$.
We note that the Coherence lemma (2.4.1 of \cite{Fra}, p.\ 50) used in the proof of 10.9.2(2) works if, in particular, the language is finite and
${\mathcal I}=[X]^{<\o}$ for some set $X$.
\kdok
\begin{fac}\label{T8089}
If $\Y$ is an infinite almost chainable $L$-structure and $F\in [Y]^{n}$ is the kernel of $\Y$, then $\f _\Y (m)\leq 2^n$, for each positive integer $m$.
\end{fac}
\dok
Let $\Y$ be $(F,<)$-chainable, where $<\in LO _{Y\setminus F}$. For $m\in \N$ we prove
\begin{equation}\label{EQ8088}
\forall K,H \in [Y]^m \;\;(K\cap F =H\cap F \Rightarrow \BK \cong \BH ),
\end{equation}
where $\BK$ and $\BH$ are the substructures of $\Y$ corresponding to $K$ and $H$ respectively.
If $K,H \in [Y]^m$ and $K\cap F =H\cap F$, then, since $|K\setminus F|=|H\setminus F|$, there is $\f \in \Pa (\la Y\setminus F ,<\ra)$
such that $\f [K\setminus F]=H\setminus F$ and by (\ref{EQ8085}), $f:=\id _F \cup \f \in \Pa (\Y )$.
Clearly we have $f[K]=H$, which implies $\BK \cong \BH$ and (\ref{EQ8088}) is true.
Now, by (\ref{EQ8088}) we have $|\{ \BK :K \in [Y]^m \}/\!\cong|\leq |P(F)|=2^n$.
\kdok
\begin{fac}\label{T8092}
Let $\Y $ and $\Z$ be $L$-structures.
If $\Y$ is almost chainable and $\Age (\Z )\subset \Age (\Y )$, then $\Z$ is almost chainable and $|\Ker (\Z )|\leq |\Ker (\Y )|$.
\end{fac}
\dok
For $|L|=1$, this is Lemma 10.9.6 of \cite{Fra}, p.\ 297, which has a straightforward generalization for arbitrary relational language $L$.
We note that 10.1.4 of \cite{Fra}, p.\ 275, which is used in the proof 10.9.6 holds for (in the notation of \cite{Fra})  $R$ and $R'$
of arbitrary signature and for $S'$ of finite signature.
\kdok
If $\,\Y\in \Mod _L (Y)$ is an infinite $(F,<)$-chainable structure, then the set
\begin{equation}\label{EQ055}
\L _\Y^F  :=  \Big\{ \la Y\setminus F, \vartriangleleft \ra : \;\vartriangleleft \in LO_{Y\setminus F}
\mbox{ and }\Y \mbox{ is } (F,\vartriangleleft )\mbox{-chainable} \Big\}
\end{equation}
is a non-empty set of linear orders
and it is easy to see that $\la Y\setminus F, \vartriangleleft \ra\in \L _\Y^F$ iff $\la Y\setminus F, \vartriangleleft \ra ^* \in \L _\Y^F$.
Theorem 9 of  \cite{Gib} gives the following description of the set $\L _\Y^F$.
\begin{te}[Gibson, Pouzet and Woodrow]\label{T8110}
If $\,\Y\in \Mod _L (Y)$ is an infinite $(F,<)$-chainable $L$-structure and $\BL :=\la Y\setminus F,< \ra$, then one of the following holds\\[-6mm]
\begin{itemize}\itemsep=-1mm
\item[{\sc (i)}] $\L _\Y^F=LO _{Y\setminus F}$, that is, each linear order $\vartriangleleft$ on $Y\setminus F$ chains $\Y$ over $F$,
\item[{\sc (ii)}] $\L _\Y ^F=\bigcup _{\BL =\BI +\BF}\Big\{ \BF + \BI ,\, \BI ^* +\BF ^*\Big\}$,
\item[{\sc (iii)}] There are finite subsets $K$ and $H$ of $\,Y\setminus F$ such that $\BL =\BK +\BM +\BH\,$ and\\[2mm]
                   $ \L _\Y^F =\bigcup _{\vartriangleleft _K \in LO_K \atop \vartriangleleft _H \in LO_H}
                   \Big\{ \la K ,\vartriangleleft _K\ra + \BM +\la H ,\vartriangleleft _H \ra,
                   \la H ,\vartriangleleft _H\ra ^*+ \BM ^* +\la K ,\vartriangleleft _K \ra^*\Big\} $.
\end{itemize}
\end{te}
\section{Almost chainable theories}\label{S3}
A complete theory $\CT \subset \Sent _L$ will be called {\it almost chainable} iff each model $\Y$ of $\CT$
is almost chainable and this notion is established by the following theorem.
\begin{te}\label{T027}
If $\CT$ is a complete $L$-theory with infinite models, then the following conditions are equivalent:

(a) All models of $\;\CT$ are almost chainable,

(b) $\CT$ has an almost chainable model,

(c) $\CT$ has a countable almost chainable model.

\noindent
If (a) is true, then there is $n\in \o$ such that $|\Ker (\Y)|=n$, for each model $\Y$ of $\CT$.
\end{te}
A proof of the theorem is given at the end of the section.
\begin{cla}\label{T029}
Let $\CK $ be a finite family of non-isomorphic $L$-structures of size $n\in \N$. Then we have

(a) For each finite set $J\subset I$ there is an $L_J$-sentence $\p _{\CK ,J}$ such that for each $\Y \in \Mod _L$ we have:
$\Y \models \p _{\CK ,J}$ iff $\;\{ \BH |L_J :H\in [Y]^n\}/\!\cong\; =\{ [ \BK |L_J ]:\BK \in \CK\} $;

(b) For the first-order theory $\CT _\CK :=\{ \p _{\CK ,J} :J\in [I]^{<\o}\}$ and each $\Y \in \Mod _L$ we have:
$\Y \models \CT _\CK$  iff  $\Age _n (\Y ) = \{ [ \BK ]:\BK \in \CK\} $.
\end{cla}
\dok
First, without loss of generality we can assume that the domain of each structure $\BK \in \CK$ is the same set, say $K$.
Let $K=\{ x_0 ,\dots ,x_{n-1}\}$ be an enumeration of its elements and $\bar x:=\la x_0 ,\dots ,x_{n-1}\ra$.

(a) For a structure $\BK \in \CK$,  let $\a_{\BK ,J} (\bar v) \!:= \bigwedge \{ \eta \in \Lit _{L_J}(\bar v): \BK \models \eta [\bar x]\} $, where
$\Lit _{L_J}(\bar v)$ is the set of all literals of $L_J$ with variables in the set $\{ v_0 ,\dots ,v_{n-1}\}$.
Then for $\Y \in \Mod _L$, $\bar y \in Y^n$ and $H:=\{ y_0 ,\dots ,y_{n-1}\}$, we have
$\Y \models \a_{\BK ,J} [\bar y]$ iff
$\{ \la x_k ,y_k\ra :k<n\}$ is an isomorphism from $\BK |L_J$ onto $\BH |L_J$.
If $\pi\in \Sym (n)$
and $\a_{\BK ,J} ^\pi (\bar v)$ is the formula obtained from $\a_{\BK ,J}$
by replacement of $v_k$ by $v_{\pi (k)}$, for all $k<n$,
then $\Y \models \a_{\BK ,J} ^\pi [\bar y]$
iff $\Y \models \a_{\BK ,J} [y_{\pi (0)}, \dots ,y_{\pi (n-1)} ]$
iff $p_\pi:=\{ \la x_k ,y_{\pi (k)}\ra :k<n\}$ is an isomorphism from $\BK |L_J$ onto $\BH |L_J$.
So, for the formula $\f_{\BK ,J}  (\bar v) := \bigvee _{\pi \in \Sym (n) }\a_{\BK ,J} ^\pi(\bar v)$
we have $\Y \models \f _{\BK ,J}[\bar y]$ iff $\BH |L_J \cong \BK |L_J$, and the equivalence in (a) is true for the formula
\begin{equation}\label{EQ051}
\textstyle
\p _{\CK ,J} := \bigwedge _{\BK \in \CK} \exists \bar v \;\f_{\BK ,J} (\bar v) \land
                \forall \bar v \; \Big( (\bigwedge _{k<l<n} v_k \neq v_l ) \Rightarrow \bigvee _{\BK \in \CK}\f_{\BK ,J} (\bar v)\Big) .
\end{equation}

(b) Let $\Y \models \CT _\CK$. Suppose that $H=\{ y_0 ,\dots ,y_{n-1}\}\in [Y]^n$ and that $\BH \not \cong \BK$, for all $\BK \in \CK$.
Then for each $\BK \in \CK$ and each $\pi \in \Sym (n)$ we have
$p_{\BK ,\pi} :=\{ \la x_k ,y_{\pi (k)}\ra :k<n\}\not\in \Iso (\BK ,\BH )$
and, since $p_{\BK ,\pi} :K\rightarrow H$ is a bijection,
there is $i_{\BK ,\pi}\in I$ such that $p_{\BK ,\pi} \not\in \Iso (\la K, R_{i_{\BK ,\pi}}^{\BK}\ra ,\la H , R_{i_{\BK ,\pi}}^{\BH }\ra )$.

Since $J:=\{ i_{\BK ,\pi} :  \BK \in \CK \land \pi \in \Sym (n)\}\in [I]^{<\omega}$ and $\Y \models \p _{\CK ,J}$,
by (a) there are  $\BK _0 \in \CK$ and $\pi _0\in \Sym (n)$
such that $p_{\BK _0,\pi _0} \in \Iso (\BK _0|L_J ,\BH |L_J )$,
which implies that $p_{\BK _0,\pi _0}\in \Iso (\la K_0, R_{i_{\BK _0,\pi _0}}^{\BK _0}\ra ,\la H , R_{i_{\BK _0,\pi _0}}^{\BH }\ra )$
and we have a contradiction. So we have proved
\begin{equation}\label{EQ048}
\forall H\in [Y]^n \;\;\exists \BK \in \CK\;\; \BH \cong \BK ,
\end{equation}
that is, $\Age _n (\Y ) \subset \{ [ \BK ]:\BK \in \CK\}$. Concerning the inclusion ``$\supset$", suppose that for some $\BK _0 \in \CK$
\begin{equation}\label{EQ047}
\forall H\in [Y]^n \;\;\BH \not\cong \BK _0
\end{equation}
and let $\CK =\{ \BK _0 ,\dots ,\BK _{s-1}\}$ be an enumeration. Then,  for each $0<r<s$ and $\pi \in \Sym (n)$, since $\BK _0 \not\cong \BK _r$, we have
$p_\pi := \{\la x_k,x_{\pi (k)}\ra: k<n \}\not\in \Iso (\BK _0 ,\BK _r)$
and, hence, there is $i_{r, \pi} \in I$ such that
\begin{equation}\label{EQ049}
 p_\pi \not\in \Iso (\la  K, R_{i_{r, \pi}}^{\BK _0}\ra ,\la  K, R_{i_{r, \pi}}^{\BK _r}\ra).
\end{equation}
Now $J:= \{i_{r, \pi}: 0<r<s \land  \pi \in \Sym (n)\}\in [I]^{<\o }$ and, since $\Y \models \p _{\CK ,J}$,
there is $H\in [Y]^n$ such that $\BH |L_J \cong \BK _0 |L_J $. By (\ref{EQ048}) and (\ref{EQ047})
there is $r>0$ such that $\BH \cong \BK _r$, which implies $\BH |L_J\cong \BK _r |L_J$ and, hence, $\BK _0 |L_J\cong \BK _r |L_J$.
Thus, there is $\pi \in \Sym (n)$ such that $p_\pi \in \Iso (\BK _0 |L_J,\BK _r|L_J)$, which in particular gives
$ p_\pi \in \Iso (\la  K, R_{i_{r, \pi}}^{\BK _0}\ra ,\la  K, R_{i_{r, \pi}}^{\BK _r}\ra)$, but this contradicts (\ref{EQ049}). So
for each $\BK \in \CK$ there is $H\in [Y]^n$ such that $\BH \cong \BK$; that is,  $\{ [ \BK ]:\BK \in \CK\}\subset \Age _n (\Y )$.

Conversely, if $\Age _n (\Y ) = \{ [ \BK ]:\BK \in \CK\}$ and $J\in [I]^{<\o}$, then for each $H\in [Y]^n$ there is $\BK \in \CK$ such that $\BH \cong \BK$;
so $\BH |L_J\cong \BK |L_J$ and, hence, $[\BH |L_J]\cong [\BK |L_J]$. In addition,  for each $\BK \in \CK$ there is $H\in [Y]^n$ such that $\BH \cong \BK$
so $\BH |L_J\cong \BK |L_J$ again and, by (a),  $\Y \models \p _{\CK ,J}$. Thus we have $\Y \models \CT _\CK$.
\kdok
\begin{cla}\label{T028}
Let $\Y $ be an infinite $L$-structure and $|\Age _n (\Y )|<\o$, for each $n\in \N$.
Then for the theory $\CT _{\rm \Age (\Y )}:=\bigcup _{n\in \N}\CT _{\Age _n(\Y )}$ and any $L$-structure $\Z$ we have

(a) $\CT _{\rm \Age (\Y )}\subset \Th (\Y )$;

(b)  $\Z \models \CT _{\rm \Age (\Y )}\,$ iff $\;\Age (\Z )=\Age (\Y )$;

(c) If $\,\Z \models \Th (\Y )$, then $\Age (\Z )=\Age (\Y )$.
\end{cla}
\dok
(a) By Claim \ref{T029}(b), for $n\!\in \!\N$ we have $\Y \!\models \CT _{\Age _n(\Y )}$
so $\CT _{\Age _n(\Y )}\subset \Th (\Y )$.

(b) $\Z \models \CT _{\rm \Age (\Y )}$
iff $\Z \models \CT _{\rm \Age ^n(\Y )}$, for all $n\in \N$;
iff (by  Claim  \ref{T029}(b))  $\Age _n(\Z )=\Age _n (\Y )$, for all $n\in \N$;
iff $\Age (\Z )=\Age (\Y )$.

(c) If $\,\Z \models \Th (\Y )$, then by (a) $\Z\models \CT _{\rm \Age (\Y )}$ and by (b) $\Age (\Z )=\Age (\Y )$.
\hfill $\Box$
\begin{cla}\label{T035}
If $\Y$ is an infinite almost chainable $L$-structure and $\,\Z \models \Th (\Y )$,
then $\Age (\Z )=\Age (\Y )$, the structure $\Z $ is almost chainable and $|\Ker (\Z )|=|\Ker (\Y )|$.
\end{cla}
\dok
By Fact \ref{T8088} we have $\Ker (\Y ) \in [Y]^n$, for some $n\in \N$,
and, by Fact \ref{T8089}, $|\Age _m (\Y )|\leq 2^n$, for all $m\in \N$. So, by Claim \ref{T028}(c) we have $\Age (\Z )=\Age (\Y )$
and, by Fact \ref{T8092}, the structure $\Z $ is almost chainable and $|\Ker (\Z )|=|\Ker (\Y )|$.
\hfill $\Box$
\begin{cla}\label{T032}
If $\CT$ is a complete almost chainable $L$-theory with infinite models and $|I|>\o$,
then $\CT$ has a countable model
and there are a countable language $L_J \subset L$ and a complete almost chainable $L _J$-theory $\CT _J$
such that
\begin{equation}\label{EQ037}
\Big|\Mod _L ^\CT (\o )/\cong \Big| =\Big|\Mod ^{\CT _J}_{L_J}(\o)/\cong \Big|.
\end{equation}
\end{cla}
\dok
Let $\Y =\la Y, \la R_i^\Y :i\in I\ra\ra\in \Mod ^\CT _L$. By Fact \ref{T8081},
there are a finite set $F=\{ a_0,\dots ,a_{n-1}\}\subset Y$, a linear order $\lhd \LO _Y$ and an $L_n$-structure $\X$
satisfying (L1)--(L3) and for each $i\in I$ there is a quantifier-free formula $\f _i (v_0 ,\dots ,v_{n_i -1})$ such that
\begin{equation}\label{EQ062}
\forall \bar y \in Y^{n_i} \;\;\Big(\bar y\in R_i^\Y \Leftrightarrow \X  \models \f _i [\bar y ]\Big) .
\end{equation}
Since there are countably many $L_n$-formulas,
there  is a partition $I=\bigcup _{j\in J}I_j$, where $|J|\leq \o$, such that, picking $i_j\in I_j$, for all $j\in J$,
we have $R_i^\Y = R_{i_j}^\Y$, for all $i\in I_j$.
So, for the $L$-sentences $\eta _{i,j}:=\forall \bar v \; (R_i (\bar v)  \Leftrightarrow R_{i_j} (\bar v))$, where $j\in J$ and $i\in I_j$, we have
$\CT _\eta := \bigcup _{j\in J}\{ \eta _{i,j} : i\in I_j \} \subset \Th _L (\Y )=\CT$.
Now, $L_J :=\la R_{i_j}:j\in J\ra \subset L$ and, using recursion,
to each $L$-formula $\f $ we adjoin an $L_J$-formula $\f _J$
in the following way:
$(v_k=v_l)_J :=v_k=v_l$;
$(R_i (v_{k_0}, \dots , v_{k_{n_i-1}}))_J := R_{i_j} (v_{k_0}, \dots , v_{k_{n_i -1}})$, for all $i\in I_j$;
$(\neg \f )_J :=\neg \f _J $; $(\f \land \p )_J :=\f _J \land \p _J$ and $(\forall v\; \f )_J :=\forall v\; \f _J $. A simple induction proves that
\begin{equation}\label{EQ036}
\forall \Z \in \Mod _L ^{\CT _\eta} \;\;\forall \f (\bar v)\in \Form _L \forall \bar z \in Z\;\;
\Big(\Z \models \f[\bar z] \Leftrightarrow \Z |L_J\models \f _J[\bar z] \Big) .
\end{equation}
We prove that, in addition, for each $\Z _1 ,\Z _2 \in \Mod _L^{\CT _\eta}$ we have
\begin{equation}\label{EQ038}
\Z _1 \cong \Z _2 \Leftrightarrow \Z _1 |L_J\cong   \Z _2 |L_J \;\;\mbox{ and }\;\;\Z _1 \equiv _L \Z _2 \Leftrightarrow \Z _1 |L_J\equiv _{L_J}   \Z _2 |L_J .
\end{equation}
The first claim is true since $\Iso (\Z _1 ,\Z _2 )=\Iso (\Z _1 |L_J, \Z _2 |L_J )$.
For the second, suppose that $\Z _1 \equiv _L \Z _2$ and $\Z _1 |L_J \models \p$, where $\p \in \Sent _{L_J}$.
Then $\p \in \Sent _L$ and $\Z _1 \models \p$, which gives $\Z _2 \models \p$
so, by (\ref{EQ036}), $\Z _2 |L_J \models \p$, because $\p _J=\p$.
Conversely, suppose that $\Z _1 |L_J\equiv _{L_J} \Z _2 |L_J$ and $\Z _1\models \f$, where $\f \in \Sent _L$.
Then, by (\ref{EQ036}), $\Z _1 |L_J \models \f _J$ and, hence, $\Z _2 |L_J \models \f _J$ so, by (\ref{EQ036}),  $\Z _2\models \f$.

Let $\CT _J :=\Th _{L_J}(\Y |L_J)$. If $\Z \in \Mod _L^{\CT }$, that is, $\Z \equiv _L \Y$,
then by (\ref{EQ038}) we have   $\Z|L_J  \equiv _{L_J} \Y|L_J $, which means that $\Z|L_J \in \Mod _{L_J}^{\CT _J }$.
So we obtain the mapping $\Lambda : \Mod _L^{\CT } \rightarrow \Mod _{L_J}^{\CT _J }$, where $\Lambda (\Z )= \Z |L_J$, for all $\Z \in \Mod _L^{\CT }$,
which is an injection, because $\CT _\eta \subset\CT $.

If $\A \in \Mod _{L_J}^{\CT _J }$, then $\Z =\la A, \la R_i^{\Z } :i\in I\ra\ra\in \Mod _L^{\CT _\eta}$,
where $R_i^{\Z }=R_{i_j}^{\A }$, for $j\in J$ and $i\in I_j$. Now $\Z |L_J =\A \equiv _{L_J}\Y |L_J$ and, by (\ref{EQ038}), $\Z \equiv _{L}\Y $,
that is, $\Z\in \Mod _L^{\CT }$ and $\Lambda $ is a surjection.
Since the mapping $\Lambda$ preserves cardinalities of structures, we have
$\Lambda [\Mod ^\CT _L (\o)]=\Mod ^{\CT _J} _{L_J}(\o)$.
By the L\"{o}wenheim-Skolem theorem there is $\A \in \Mod ^{\CT _J}_{L_J}(\o)$ and
$\Lambda ^{-1}(\A )$ is a countable model of $\CT$.

By (\ref{EQ038}), the mapping $\Lambda$ preserves the isomorphism relation and (\ref{EQ037}) is true.
By (\ref{EQ062}) the reduct $\Y|L_J$ is simply definable in $\X$
and, by Fact \ref{T8081}, it is almost chainable.
By  Claim \ref{T035}, the theory $\CT _J=\Th _{L_J}(\Y |L_J)$ is almost chainable.
\kdok
\noindent
{\bf Proof of Theorem \ref{T027}.}
The implication (a) $\Rightarrow$ (c) follows from  Claim \ref{T032}, the implication (c) $\Rightarrow$ (b) is trivial
and (b) $\Rightarrow$ (a) follows from  Claim \ref{T035}.
\hfill $\Box$
\section{Vaught's Conjecture}\label{S4}
In this section we confirm Vaught's Conjecture for almost chainable theories. More precisely,
the whole section is devoted to a proof of the following statement.
\begin{te}\label{T030}
If $\CT$ is a complete almost chainable theory having infinite models, then
$I(\CT ,\o)\in \{ 1 , \c \}$. In addition, the theory $\CT$ is $\o$-categorical iff it has a countable model
which is chained by an $\o$-categorical linear order over its kernel.
\end{te}
So, let $\CT$ be a complete almost chainable $L$-theory having infinite models.
By Theorem \ref{T027}, there is $n\in \o$ such that each model of $\CT$ has the kernel of size $n$ and,
by  Claim \ref{T032}, w.l.o.g.\ we suppose that $|L|\leq \o$, which gives $\Mod_L^{\CT} (\o)\neq \emptyset$.
As above, let $L_n$ denote the language $\la R, U_0 ,\dots U_{n-1}\ra$, where $R$ is a binary and $U_j$'s are unary symbols.
In the sequel, for $\Y \in \Mod_L (\o)$, by $[\Y ]$ we denote the set $\{ \Y ' \in \Mod_L (\o): \Y '\cong \Y\}$ and similarly
for the structures from $\Mod_{L_n} (\o)$.

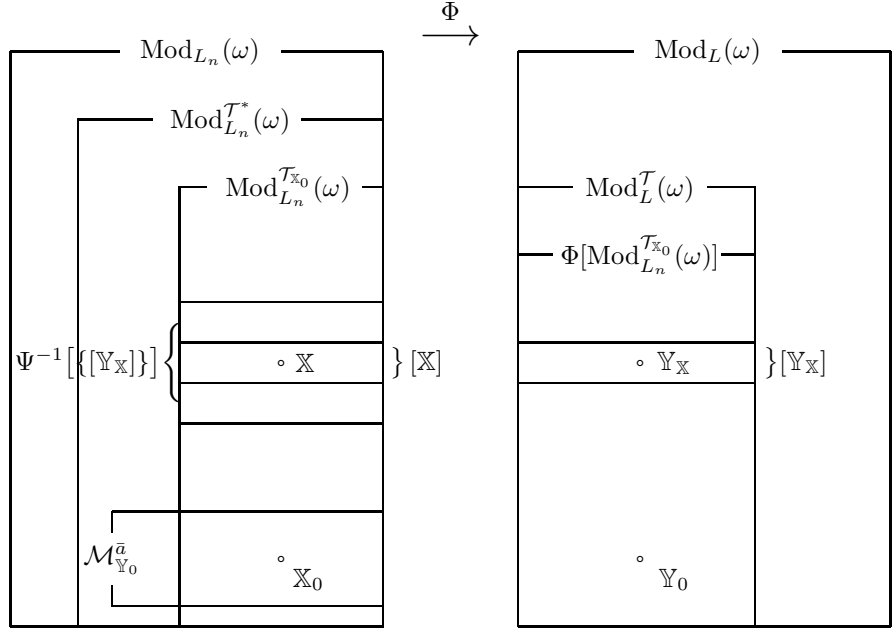
\begin{figure}[h]
\begin{center}
\unitlength 0.9mm 
\linethickness{0.5pt}
\ifx\plotpoint\undefined\newsavebox{\plotpoint}\fi 
\begin{picture}(140,100)(0,0)
\put(5,5){\line(1,0){55}}
\put(80,5){\line(1,0){55}}
\put(20,8){\line(1,0){40}}
\put(20,22){\line(1,0){40}}
\put(30,35){\line(1,0){30}}
\put(30,41){\line(1,0){30}}
\put(80,41){\line(1,0){35}}
\put(30,47){\line(1,0){30}}
\put(80,47){\line(1,0){35}}
\put(30,53){\line(1,0){30}}
\put(80,60){\line(1,0){5}}
\put(115,60){\line(-1,0){5}}
\put(30,70){\line(1,0){4}}
\put(57,70){\line(1,0){3}}
\put(80,70){\line(1,0){8}}
\put(115,70){\line(-1,0){7}}
\put(15,80){\line(1,0){11}}
\put(48,80){\line(1,0){12}}
\put(5,90){\line(1,0){16}}
\put(45,90){\line(1,0){15}}
\put(80,90){\line(1,0){17}}
\put(119,90){\line(1,0){16}}
\put(5,5){\line(0,1){85}}
\put(15,5){\line(0,1){75}}
\put(30,5){\line(0,1){65}}
\put(20,19){\line(0,1){3}}
\put(20,8){\line(0,1){3}}
\put(60,5){\line(0,1){85}}
\put(80,5){\line(0,1){85}}
\put(115,5){\line(0,1){65}}
\put(135,5){\line(0,1){85}}
\put(45,15){\circle{0.8}}
\put(98,15){\circle{0.8}}
\put(45,44){\circle{0.8}}
\put(98,44){\circle{0.8}}
\small
\put(20,15){\makebox(0,0)[cc]{$\CM _{\Y _0}^{\bar a} $}}
\put(49,12){\makebox(0,0)[cc]{$\X_0 $}}
\put(103,12){\makebox(0,0)[cc]{$\Y_0 $}}
\put(18,44){\makebox(0,0)[cc]{$ \Psi ^{-1}\big[ \{[\Y _\X]\}\big] \Bigg\{$}}
\put(48.5,44){\makebox(0,0)[cc]{$\X $}}
\put(65,44){\makebox(0,0)[cc]{$\big\}\, [\X] $}}
\put(103,44){\makebox(0,0)[cc]{$\Y _\X $}}
\put(121,44){\makebox(0,0)[cc]{$\big\} [\Y _\X] $}}
\put(98,60){\makebox(0,0)[cc]{$ \Phi [ \Mod _{L_n}^{\CT _{\X _0}}(\o )]$}}
\put(46,70){\makebox(0,0)[cc]{$\Mod _{L_n}^{\CT _{\X _0}}(\o ) $}}
\put(98,70){\makebox(0,0)[cc]{$\Mod _{L}^{\CT }(\o ) $}}
\put(37,80){\makebox(0,0)[cc]{ $\Mod _{L_n}^{\CT ^*}(\o )$}}
\put(33,90){\makebox(0,0)[cc]{$ \Mod _{L_n}(\o ) $}}
\put(108,90){\makebox(0,0)[cc]{$\Mod _{L}(\o ) $}}
\put(70,96){\makebox(0,0)[cc]{$\Phi $}}
\Large
\put(70,92){\makebox(0,0)[cc]{$\longrightarrow$}}
\end{picture}
\end{center}
\caption{The mappings $\Phi$ and $\Psi$}\label{F000}
\end{figure}

Following the architecture of the proof of the corresponding statement from \cite{K}
we divide the proof into two subsections.
In ``Preliminaries" we take an arbitrary countable model $\Y _0$ of $\CT$
and a linear order with $n$ unary predicates $\X _0$ such that $\Pa (\X _0)\subset \Pa (\Y _0)$ (see Figure \ref{F000})
and describe the cardinal argument which will be used in our proof.
In ``Proof", distinguishing some cases, taking convenient structures $\Y _0$ and $\X _0$ and using that cardinal argument, we prove Theorem \ref{T030}.

\subsection{Preliminaries}
For convenience,
let $\Delta _n:=\{ \la x_0,\dots ,x_{n-1}\ra \in \o ^n:\bigvee _{k<l<n}x_k=x_l\}$
and, for an $n$-tuple $\bar a:=\la a_0,\dots ,a_{n-1}\ra\in \o ^n$, let us define $F_{\bar a}:=\{ a_0,\dots ,a_{n-1}\}$.

We fix a model $\Y _0=\la \o , \la R _i ^{\Y _0}:i\in I \ra \ra\in \Mod _{L}^{\CT }(\o )$ and an enumeration of its kernel, $\Ker (\Y _0)=\{ a_0,\dots ,a_{n-1}\} $.
By Fact \ref{T8081} there is a linear order $\prec \in LO _\o$ such that, defining $\bar a =\la a_0,\dots ,a_{n-1}\ra$
and $\X _0:= \la \o ,\prec , \{a_0\},\dots ,\{a_{n-1}\}\ra$,
$$
\la \o ,\prec \ra =\{a_0\}+ \dots +\{a_{n-1}\}+(\o \setminus F_{\bar a}) \;\;\mbox{ and } \;\;\Pa (\X _0)\subset \Pa (\Y _0).
$$
Thus, the structure $\Y _0$ is $(F_{\bar a},\prec\upharpoonright (\o \setminus F_{\bar a}))$-chainable.
Let $\CT _{\X _0}$ denote the complete theory of $\X_0$, $\Th _{L_n}(\X _0)$.
The structure $\X_0$ has the following  properties expressible by first order
sentences of the language $L_n$:\\[-7mm]
\begin{itemize}\itemsep -1mm
\item[(i)] The interpretation of $R$ is a linear order,
\item[(ii)] The interpretations of the relations $U_k$, $k<n$, are different singletons,
\item[(iii)] These singletons are ordered as the indices of $U_k$'s
(that is, the $L_n$-sentence $ \bigwedge _{k<l<n} \forall u,v \; (U_k(u) \land U_l(v)\Rightarrow R(u,v))$ is true in $\X _0$),
\item[(iv)] The union of these singletons is an initial segment of the linear order; that is
$\X _0 \models \forall u,v \; ((U_{n-1}(u)\land \bigwedge _{k<n}\neg U_k(v))\Rightarrow R(u,v))$.
\end{itemize}
So, if $\CT ^*$ is the set of the $L_n$-sentences expressing (i)--(iv), then $\CT ^*\subset \CT _{\X _0}$ and
\begin{eqnarray*}
\X _0 \in \CM ^{\bar a}_{\Y _0} & := & \Big\{ \la \o ,\vartriangleleft , \{a_0\},\dots ,\{a_{n-1}\}\ra \in \Mod _{L_n}^{\CT ^*}(\o ):\\
                       &    & \;\;\Y _0 \mbox{ is } (F_{\bar a}, \vartriangleleft \upharpoonright (\o \setminus F_{\bar a})) \mbox{-chainable}\Big\}.
\end{eqnarray*}
By Fact \ref{T8081} the structure $\Y _0$ is simply definable in the $L_n$-structure $\X _0$.
Thus, for each $i\in I$ there is a quantifier free $L_n$-formula $\f _i (v_0 ,\dots ,v_{n_i-1})$ such that
\begin{equation}\label{EQ024}
\forall \bar x\in \o ^{n_i} \;\; \Big(\bar x\in R _i ^{\Y _0}\Leftrightarrow  \X _0 \models \f _i[\bar x]\Big).
\end{equation}
Generally speaking, using the $L_n$-formulas $\f_i$, $i\in I$,
to each $L_n$-structure $\X\in \Mod _{L_n}(\o )$ we can adjoin the $L$-structure
$\Y _\X:=\la \o , \la R _i ^{\Y _\X}:i\in I \ra \ra \in \Mod _{L}(\o )$,
where, for each $i\in I$, the relation $R _i ^{\Y _\X}$ is defined in the structure $\X$
by the formula $\f_i$, that is,
\begin{equation}\label{EQ042}
\forall \bar x\in \o ^{n_i} \;\; \Big(\bar x\in R _i ^{\Y _\X }\Leftrightarrow  \X  \models \f _i[\bar x]\Big).
\end{equation}
\begin{cla}\label{T034}
For each structure $\Y _0\in \Mod _{L}^{\CT }(\o )$, each enumeration $\Ker (\Y _0)=\{ a_0,\dots ,a_{n-1}\}$
each
structure $\X_0 \in \CM _{\Y _0}^{\bar a}$ and each choice of formulas $\f _i$, $i\in I$, satisfying (\ref{EQ024}),
defining $\Y _\X$ by (\ref{EQ042}), for $\X\in \Mod _{L_n}(\o )$,  we have\\[-5.5mm]
\begin{itemize}\itemsep -1mm
\item[\rm (a)] The mapping
$\Phi : \Mod _{L_n}(\o )\rightarrow \Mod _{L}(\o )$, defined by
$
\Phi (\X )=\Y _\X,
$
for each $\X \in \Mod _{L_n}(\o )$,
preserves elementary equivalence and isomorphism;
moreover, $\Iso (\X _1, \X _2)\subset \Iso (\Y _{\X _1}, \Y _{\X _2})$, for all  $\X _1, \X _2 \in \Mod _{L_n}(\o )$;
\item[\rm (b)] The mapping $\Psi : \Mod _{L_n}^{\CT _{\X _0}}(\o )/\!\cong  \;\rightarrow \Mod _{L}^{\CT }(\o )/\! \cong$, given by
$
\Psi ([\X ])=[\Y _\X],
$
for all $[\X ]\in \Mod _{L_n}^{\CT _{\X _0}}(\o )/\!\cong $,
is well defined.
\end{itemize}
\end{cla}
\dok
(a) By recursion on the construction of $L$-formulas
to each $L$-formula $\f (\bar v)$ we adjoin an $L_n$-formula $\f ^*(\bar v)$ in the following way:
$(v_k=v_l)^* :=v_k=v_l$, $R_i (v_{k_0}, \dots , v_{k_{n_i-1}})^*:= \f _i (v_{k_0}, \dots , v_{k_{n_i-1}})$ (replacement of $v_j$ by $v_{k_j}$ in $\f_i$),
$(\neg \f)^* := \neg \f^*$, $(\f\land \p)^* := \f^*\land \p^*$ and $(\exists v_k\; \f)^* := \exists v_k\; \f^*$. A routine induction shows that, writing
$\bar v$ instead of  $v_0,\dots ,v_{n-1}$, we have (see \cite{Hodg}, p.\ 216)
\begin{equation}\label{EQ025}
\forall \X \in \Mod _{L_n}(\o ) \;\;\forall \f (\bar v)\in \Form _L\;\forall \bar x \in \o ^n \;
\Big( \X \models \f ^*[\bar x]\Leftrightarrow \Y_\X \models \f[\bar x]\Big).
\end{equation}
Let $\X _1, \X _2 \in \Mod _{L_n}(\o )$. If $\X_1 \equiv \X_2$, then for an $L$-sentence $\f$ we have:
$\Y _{\X _1} \models \f$
iff  $\X _1\models \f ^*$ (by (\ref{EQ025}))
iff $\X _2 \models \f ^*$ (since $\X _1 \equiv \X_2$)
iff $\Y _{\X _2} \models \f$ (by (\ref{EQ025}) again). So, $\Y _{\X _1} \equiv \Y _{\X_2} $ and the mapping $\Phi $
preserves elementary equivalence.

If $f:\X_1 \rightarrow  \X _2$ is an isomorphism,
then by (\ref{EQ042}) and since  isomorphisms preserve all formulas in both directions, for each  $i\in I$ and $\bar x\in \o ^{n_i}$ we have:
$\bar x\in R_i ^{\Y_{\X _1}}$
iff $\X _1 \models \f _i [\bar x]$
iff $\X _2 \models \f _i [f\bar x]$
iff $f \bar x\in R_i ^{\Y_{\X _2}}$.
Thus $f\in \Iso (\Y _{\X_1}, \Y _{ \X _2})$.

(b) For $\X \in \Mod _{L_n}^{\CT _{\X _0}}(\o )$ we have $\X \equiv \X_0$, which, by (a),  (\ref{EQ024}) and (\ref{EQ042}),
implies that $\Phi (\X )=\Y _\X \equiv \Y _{\X_0}=\Y _0$.
So, since $\Y_0\models \CT$, we have $\Phi (\X )\in \Mod _{L}^{\CT }(\o )$ and, thus,
\begin{equation}\label{EQ043}
\Phi [ \Mod _{L_n}^{\CT _{\X _0}}(\o )] \subset \Mod _{L}^{\CT }(\o ).
\end{equation}
Assuming that $\X_1, \X _2 \in \Mod _{L_n}^{\CT _{\X _0}}(\o )$ and $\X _1\cong \X _2$, by (a) we have $\Y_{\X _1}\cong \Y_{\X _2}$,
that is $[\Y_{\X _1}]=[\Y_{\X _2}]$. So, the mapping $\Psi $ is well defined.
\kdok
\noindent
Thus, by Claim \ref{T034}(b), if $I(\CT _{\X _0},\o )=\c$, then for a proof that $I(\CT ,\o)=\c$
it is sufficient to show that the mapping $\Psi$ is at-most-countable-to-one,
which will be true if for each $\X\in \Mod _{L_n}^{\CT _{\X _0}}(\o )$ we have $|\Psi ^{-1}[ \{[\Y _\X]\}]\leq \o$.
We note that, by Example 4.2 of \cite{K}, it is possible that $|\Psi ^{-1}[ \{[\Y _\X]\}]= \o$.

Now, let $\X\in \Mod _{L_n}^{\CT _{\X _0}}(\o )$.
Then we have $\X\models \CT ^*$
and, hence, there is $\bar b:=\la b_0,\dots ,b_{n-1}\ra\in \o ^n \setminus \Delta _n$
such that $\X= \la \o ,\prec _\X, \{b_0\},\dots ,\{b_{n-1}\}\ra$ and $\Y _\X$ is definable in $\X$ by (\ref{EQ042}).
So, by Fact \ref{T8081}, the structure $\Y _\X$ is $(F_{\bar b}, \prec _\X\upharpoonright (\o \setminus F_{\bar b}))$-chainable and
(see (\ref{EQ055}))
\begin{equation}\label{EQ054}
\BL _\X := \la \o \setminus F_{\bar b}, \prec _\X\upharpoonright (\o \setminus F_{\bar b}) \ra\in \L _{\Y _\X}^{F_{\bar b}}.
\end{equation}
For an $n$-tuple $\bar c:=\la c_0,\dots ,c_{n-1}\ra\in \o ^n\setminus \Delta _n$ let us define
\begin{eqnarray}
\CM ^{\bar c}_{\Y _\X} & := & \Big\{ \la \o ,\prec , \{c_0\},\dots ,\{c_{n-1}\}\ra \in \Mod _{L_n}^{\CT ^*}(\o ):\nonumber \\
                       &    & \;\;\Y _\X \mbox{ is } (F_{\bar c}, \prec \upharpoonright (\o \setminus F_{\bar c})) \mbox{-chainable}\Big\},\label{EQ056}\\
\CM _{\Y _\X}          & := & \textstyle \bigcup _{\bar c\in \o ^n\setminus \Delta _n}\CM ^{\bar c}_{\Y _\X}\label{EQ057}.
\end{eqnarray}
Thus, $\X \in \CM ^{\bar b}_{\Y _\X}\subset \CM _{\Y _\X}$.
For $\CM \subset \Mod _{L_n}^{\CT ^*}(\o )$, let $\CM ^{\cong }:=\{ [\A ]:\A \in \CM\}$.
\begin{cla}\label{T036}
For each structure $\X\in \Mod _{L_n}^{\CT _{\X _0}}(\o )$ we have
$$
\Psi ^{-1}\Big[ \{[\Y _\X]\}\Big] \subset \CM _{\Y_{\X  }}^{\cong}\cap \Mod _{L_n}^{\CT _{\X _0}}(\o )/\cong.
$$
\end{cla}
\dok
Let $\X\in \Mod _{L_n}^{\CT _{\X _0}}(\o )$.
Then, by (\ref{EQ043}) we have $\Y _\X \in \Mod _{L}^{\CT }(\o )$; so $\Psi ([\X])=[\Y _\X ]\in \Mod _{L}^{\CT }(\o )/\cong$
and $\Psi ^{-1}[\{ [\Y _\X]\}]\subset \dom (\Psi)=\Mod _{L_n}^{\CT _{\X _0}}(\o )/\cong $.

Let $[\X _1]\in \Psi ^{-1}[\{[\Y _\X]\}]$.
Then $[\X _1]\in \Mod _{L_n}^{\CT _{\X _0}}(\o )/\cong$
and, since the set $\Mod _{L_n}^{\CT _{\X _0}}(\o )$ is closed under $\cong$, we have $\X _1\in \Mod _{L_n}^{\CT _{\X _0}}(\o )$.
This implies that $\X _1 \models \CT ^*$
and, hence, $\X _1 =\la \o ,\prec _{\X _1}, \{c_0 \},\dots ,\{c_{n-1}\}\ra$,
for some $\bar c:=\la c_0,\dots ,c_{n-1}\ra\in \o ^n\setminus \Delta _n$.
Since $\X _1\in \Mod _{L_n}^{\CT _{\X _0}}(\o )$, by (\ref{EQ042}) for $i\in I$ we have
\begin{equation} \label{EQ053}
\forall \bar x\in \o ^{n_i}\;\;( \bar x\in R_i^{\Y _{\X _1}} \Leftrightarrow \X _1 \models \f _i [\bar x]).
\end{equation}
Also we have $[\Y _{\X _1}]=\Psi ([\X _1])=[\Y _\X] $,
so there is $f\in \Iso (\Y _{\X } , \Y _{\X _1})$ and we prove that
$[\X _1]\in \;\CM _{\Y _{\X }} ^{\cong}$.

Clearly, $\X _2:=\la \o , f^{-1}[\prec _{\X _1}],\{f^{-1}(c_0)\},\dots ,\{f^{-1}(c_{n-1})\} \ra \cong\X_1$
and $f\in \Iso (\X _2 ,\X _1)$.
For $i\in I$ and $\bar x\in \o ^{n_i}$ we have
$\bar x\in R_i^{\Y _\X}$
iff  $f\bar x\in R_i^{\Y _{\X _1}}$ (since $f\in \Iso (\Y _\X , \Y _{\X _1})$),
iff  $\X _1 \models \f _i [f\bar x]$ (by \ref{EQ053}),
iff  $\X _2\models \f _i [\bar x]$ (since $f\in \Iso (\X _2 ,\X _1 )$).
Thus $\bar x\in R_i^{\Y _\X}$ iff $\X _2\models \f _i [\bar x]$, for all $\bar x\in \o^{n_i}$,
so, by Fact \ref{T8081},
the structure $\Y _\X$ is $(F_{f^{-1}\bar c},f^{-1}[\prec _{\X _1}]\upharpoonright(\o \setminus F_{f^{-1}\bar c}) )$-chainable.
Now, $\X _2 \in  \CM _{\Y _\X} ^{f^{-1}\bar c}$
and, hence, $[\X _1]=[\X _2]\in (\CM _{\Y _\X} ^{f^{-1}\bar c})^{\cong}\;\subset \CM _{\Y _{\X }} ^{\cong}$.
Thus $\Psi ^{-1}[ \{[\Y _\X]\}] \subset \CM _{\Y_{\X  }}^{\cong}$.
\kdok

\noindent
The following folklore statement will be used in our case analysis as well.
\begin{cla}\label{T019'}
If some structure $\Y \in \Mod _{L}^{\CT }(\o )$ is simply definable in an $\o$-categorical structure $\X$ with domain $\o$,
then $\Y $ is an $\o$-categorical structure and $I(\CT ,\o ) = 1$.
\end{cla}
\dok
By the theorem of Engeler, Ryll-Nardzewski and Svenonius (see \cite{Hodg}, p.\ 341),
the automorphism group of $\X $ is oligomorphic; that is, for each $n\in \N$ we have
$|\o ^n /\!\sim _{\X ,n}|<\o$, where $\bar x \sim_{\X ,n} \bar y$ iff $f\bar x =\bar y$, for some $f\in \Aut (\X )$.

As in Claim \ref{T034}(a) we prove that $\Aut (\X )\subset \Aut (\Y )$, which implies that for $n\in \N$
and each $\bar x , \bar y \in \o ^n$ we have $\bar x \sim_{\X ,n} \bar y \Rightarrow \bar x \sim_{\Y ,n} \bar y$.
Thus $|\o ^n /\!\sim _{\Y ,n}| \leq |\o ^n /\!\sim _{\X ,n}|<\o$, for all $n\in \N$, and, since $|L|\leq \o$, using the same theorem we conclude that
$\Y $ is an $\o$-categorical $L$-structure.
\hfill$\Box$
\subsection{Proof}
First we prove that $|\Mod _{L}^{\CT }(\o )/ \!\cong |\in \{ 1,\c\}$,
using definitions and notation from ``Preliminaries" and distinguishing the following cases.\\[-2mm]

\noindent
{\bf Case A:} {\it There exist a structure $\Y_0 \in \Mod _{L}^{\CT }(\o )$, an enumeration of its kernel, $\Ker (\Y _0)=\{ a_0,\dots ,a_{n-1}\}$,
and a  structure $\X _0\in \CM _{\Y _0} ^{\bar a}$ such that the theory $\CT _{\X _0}$ is $\o$-categorical.}
Then by Fact \ref{T8081} the structure $\Y _0$ is simply definable in $\X _0$ and by Claim \ref{T019'} we have $I(\CT ,\o ) = 1$.

In particular, Case A appears if there is a structure $\Y \in \Mod _{L}^{\CT }(\o )$ satisfying condition {\sc (i)} of Theorem \ref{T8110}:
$\Y$ is $F$-chainable and $\L _\Y^F=LO _{\o\setminus F}$. Then,
taking an enumeration $F=\{ a_0,\dots ,a_{n-1}\}$,
the relations $R^\Y_i$ of the structure $\Y$ are definable in the
structure $\X :=\la \o , \{a_0\},\dots ,\{a_{n-1}\}\ra$ of the unary language $L':=\la U_0,\dots ,U_{n-1}\ra$
by quantifier free $L'$-formulas and,
since the structure $\X$ is $\o$-categorical, $\Y$ is $\o$-categorical as well; so, $I(\CT ,\o ) = 1$ again.
We note that such structures are called {\it finitist} by Fra\"{\i}ss\'{e}, see \cite{Fra}, p.\ 292.\\[-2mm]

\noindent
{\bf Case B:} {\it For each structure $\Y_0 \in \Mod _{L}^{\CT }(\o )$, each enumeration of its kernel $\Ker (\Y _0)=\{ a_0,\dots ,a_{n-1}\}$
and each structure $\X _0\in \CM _{\Y _0} ^{\bar a}$, the theory $\CT _{\X _0}$ is not $\o$-categorical; so, by Theorem \ref{T002},
$|\Mod _{L_n}^{\CT _{\X _0 }}(\o )/\!\cong|=\c$, for all $\X _0\in \CM _{\Y _0} ^{\bar a}$.}
Then, by the remark from Case A concerning condition {\sc (i)} of Theorem \ref{T8110}, we have
\begin{equation} \label{EQ040}
\forall \Y  \in \Mod _{L}^{\CT }(\o ) \;\; \L _\Y^{\Ker (\Y )} \neq LO _{\o \setminus \Ker (\Y )} ,
\end{equation}
and we prove that $|\Mod ^\CT _L (\o )/\!\cong|=\c$, distinguishing the following two subcases.\\[-2mm]

\noindent
{\bf Subcase B1:}
{\it There exist a structure $\Y_0 \in \Mod _{L}^{\CT }(\o )$, an enumeration of its kernel, $\Ker (\Y _0)=\{ a_0,\dots ,a_{n-1}\}$,
and a  structure $\X _0\in \CM _{\Y _0} ^{\bar a}$
such that the linear order
$\BL _{\X _0} := \la \o \setminus F_{\bar a}, \prec _{\X _0}\upharpoonright (\o \setminus F_{\bar a}) \ra\in \L _{\Y _0}^{F_{\bar a}}$ has at least one end-point.}

Then we take such $\Y_0$, $\bar a$ and $\X _0$ and notice that $\X _0 \models \CT ^*$ and that the mentioned property of $\BL _{\X _0}$ gives a first order property of $\X_0$. Namely, $\X _0\models \theta _0 \lor \theta _1$, where
\begin{eqnarray}
\theta _0 & := & \exists v \;\forall u \;\Big(U_{n-1}(u)\Rightarrow R(u,v)\land \neg \exists w \;(R(u,w) \land R(w,v))\Big),\label{EQ059}\\
\theta _1 & := & \exists v \;\forall u \; \Big(\neg u=v \Rightarrow R(u,v)\Big).\label{EQ060}
\end{eqnarray}
Now we have $|\Mod _{L_n}^{\CT _{\X _0}}(\o )/\!\cong|=\c$
and, by Claim \ref{T034}(b), for a proof that $|\Mod ^\CT _L (\o )/\!\cong|=\c$
it is sufficient to show that the mapping $\Psi$ is at-most-countable-to-one. This will follow from the following claim and
Claim \ref{T036}.
\begin{cla}\label{T023}
$\Big|\CM _{\Y_{\X  }}^{\cong}\cap \Mod _{L_n}^{\CT _{\X _0}}(\o )/\cong\Big|\leq \o$, for all $\X\in \Mod _{L_n}^{\CT _{\X _0}}(\o )$.
\end{cla}
\dok
Let $\X\!\in \!\Mod _{L_n}^{\CT _{\X _0}}(\o )$. By (\ref{EQ057}) it is sufficient to show that for each $\bar c\!\in \o ^n\setminus \Delta _n$
we have
\begin{equation} \label{EQ058}
\Big| (\CM ^{\bar c}_{\Y _\X})^{\cong } \cap \Mod _{L_n}^{\CT _{\X _0}}(\o )/\cong\Big| \leq \o.
\end{equation}
Let $\X _1= \la \o ,\prec _{\X _1}, \{c_0\},\dots ,\{c_{n-1}\}\ra\in \CM ^{\bar c}_{\Y _\X}$.
Then by (\ref{EQ056}) and (\ref{EQ055}) we have
$$
\BL _{\X _1} := \la \o \setminus F_{\bar c}, \prec _{\X _1}\upharpoonright (\o \setminus F_{\bar c}) \ra\in \L _{\Y _\X}^{F_{\bar c}}.
$$
First, if the set $\L _{\Y _\X}^{F_{\bar c}}$ satisfies condition {\sc (iii)} of Theorem \ref{T8110},
then we have $(\L _{\Y _\X}^{F_{\bar c}})^{\cong }=\{ [\BL _{\X _1}], [\BL _{\X _1}^*]\}$
(because all ``$\BK +\BM +\BH$-sums" are isomorphic
and all ``$\BH ^* +\BM ^* +\BK ^*$-sums" are isomorphic).
Thus each structure $\X _2 \in \CM ^{\bar c}_{\Y _\X}$ consists of
the initial part, $\{c_0\}+\dots +\{c_{n-1}\}$,
labeled by the unary relations $U_j^{\X _2} =\{ c_j \}$, $j<n$,
and a final part, which is either isomorphic to the linear order $\BL _{\X _1}$ or to its reverse, $\BL _{\X _1}^*$.
So we have $|(\CM ^{\bar c}_{\Y _\X})^{\cong }|\leq 2$ and (\ref{EQ058}) is true.

Otherwise, by (\ref{EQ040}) and Theorem \ref{T8110},
$\L _{\Y_\X}^{F_{\bar c}} =\bigcup _{\BL _{\X _1}=\BI +\BF}\{ \BF + \BI ,\, \BI ^* +\BF ^*\}$.

Let $\X _2= \la \o ,\prec _{\X _2}, \{c_0\},\dots ,\{c_{n-1}\}\ra\in \CM ^{\bar c}_{\Y _\X}\cap \Mod _{L_n}^{\CT _{\X _0}}(\o )$.
Then $\BL _{\X _2} := \la \o \setminus F_{\bar c}, \prec _{\X _2}\upharpoonright (\o \setminus F_{\bar c}) \ra\in \L _{\Y _\X}^{F_{\bar c}}$
and, hence, there is a cut $\{ \BI,\BF\}$ in $\BL _{\X _1}$ (i.e.\ a decomposition $\BL _{\X _1}=\BI +\BF$)
such that $\BL _{\X _2} =\BF +\BI$ or $\BL _{\X _2} =\BI ^* +\BF ^*$.
Suppose that $\BI,\BF\neq\emptyset$,
that $\BI$ does not have a largest element
and that $\BF$ does not have a smallest element.
Then  $\BF +\BI$ and $\BI ^* +\BF ^*$ are linear orders without end points.
But, since $\X _2\in \Mod _{L_n}^{\CT _{\X _0}}(\o )$ we have $\X _2 \models \theta _0 \lor \theta _1$
and, hence, the linear order $\BL _{\X _2}$ must have at least one end-point,
which gives a contradiction.

So, for each $\X _2\in \CM ^{\bar c}_{\Y _\X}\cap \Mod _{L_n}^{\CT _{\X _0}}(\o )$ we have
$\BL _{\X _2} =\BF +\BI$ or $\BL _{\X _2} =\BI ^* +\BF ^*$,
where $\BI$ has a largest element or $\BF$ has a smallest element.
Since such cuts $\{ \BI,\BF\}$ in $\BL _{\X _1}$ are defined by the elements of the set  $\o \setminus F_{\bar c}$,
there are countably many of them.
Thus $|\CM ^{\bar c}_{\Y _\X}\cap \Mod _{L_n}^{\CT _{\X _0}}(\o )|=\o $,
which implies (\ref{EQ058}), since each class from the set $(\CM _{\Y_{\X  }}^{\bar c})^{\cong}\cap \Mod _{L_n}^{\CT _{\X _0}}(\o )/\cong$
has a representative in $\CM ^{\bar c}_{\Y _\X}\cap \Mod _{L_n}^{\CT _{\X _0}}(\o )$.
\kdok

\noindent
{\bf Subcase B2:} {\it  For each structure $\Y_0 \in \Mod _{L}^{\CT }(\o )$, each enumeration of its kernel, $\Ker (\Y _0)=\{ a_0,\dots ,a_{n-1}\}$,
and each  structure $\X _0\in \CM _{\Y _0} ^{\bar a}$,
the linear order
$\BL _{\X _0} := \la \o \setminus F_{\bar a}, \prec _{\X _0}\upharpoonright (\o \setminus F_{\bar a}) \ra\in \L _{\Y _0}^{F_{\bar a}}$
is a linear order without end points.}

Then we fix arbitrary $\Y_0\in \Mod _{L}^{\CT }(\o )$ and $\X _0\in \CM _{\Y _0} ^{\bar a}$, where $\Ker (\Y _0)=F_{\bar a}$.

Again we have $|\Mod _{L_n}^{\CT _{\X _0}}(\o )/\!\cong|=\c$ and,
as in Subcase B1, the equality $|\Mod ^\CT _L (\o )/\!\cong|=\c\,$
will follow from Claims \ref{T034}(b), \ref{T036} and the next claim.
\begin{cla}\label{T037}
$\Big|\CM _{\Y_{\X  }}^{\cong}\cap \Mod _{L_n}^{\CT _{\X _0}}(\o )/\cong\Big|\leq \o$, for all $\X\in \Mod _{L_n}^{\CT _{\X _0}}(\o )$.
\end{cla}
\dok
Let $\X\!\in \!\Mod _{L_n}^{\CT _{\X _0}}(\o )$. By (\ref{EQ057}) it is sufficient to show that for each $\bar c\!\in \o ^n\setminus \Delta _n$
we have $| (\CM ^{\bar c}_{\Y _\X})^{\cong } | \leq 2$.

Let $\X _1= \la \o ,\prec _{\X _1}, \{c_0\},\dots ,\{c_{n-1}\}\ra\in \CM ^{\bar c}_{\Y _\X}$. Then, by our assumption,
$\BL _{\X _1} = \la \o \setminus F_{\bar c}, \prec _{\X _1}\upharpoonright (\o \setminus F_{\bar c}) \ra$
is a linear order without end points  and $\BL _{\X _1} \in \L _{\Y _\X}^{F_{\bar c}}$.

Suppose that the set $\L _{\Y _\X}^{F_{\bar c}}$ satisfies condition {\sc (ii)} of Theorem \ref{T8110};
that is, $\L _{\Y_\X}^{F_{\bar c}} =\bigcup _{\BL _{\X _1}=\BI +\BF}\{ \BF + \BI ,\, \BI ^* +\BF ^*\}$.
Then, taking an arbitrary $x\in \o\setminus F_{\bar c}$
we have $\BL _{\X _1} =(-\infty ,x)_{\BL _{\X _1}}+ [x, \infty)_{\BL _{\X _1} }=:\BI +\BF$;
and, since $\BL _{\X _1}$ is a linear order without end points, $\BI ,\BF\neq \emptyset$.
Let $\X _2= \la \o ,\prec _{\X _2}, \{c_0\},\dots ,\{c_{n-1}\}\ra$,
where $\prec _{\X _2}$ is the linear order on $\o$
such that $\la \o , \prec _{\X _2}\ra =\{c_0\}+\dots +\{c_{n-1}\}+ \BF +\BI$.
Then the structure $\Y _\X$ is $(F_{\bar c},\prec _{\X _2}\upharpoonright (\o \setminus F_{\bar c}) )$-chainable
and, hence, $\X _2\in \CM ^{\bar c}_{\Y _\X}$.
But  the linear order $\BL _{\X _2 } := \la \o \setminus F_{\bar c}, \prec _{\X _2 }\upharpoonright (\o \setminus F_{\bar c}) \ra=\BF +\BI$
has a smallest element,
which contradicts the assumption of Subcase B2.

Thus there are finite sets $K,H\subset \o \setminus F_{\bar c}$ such that $\BL _{\X _1} =\BK +\BM +\BH\,$ and
$$\textstyle
\L _{\Y _\X}^{F_{\bar c}} =\bigcup _{\vartriangleleft _K \in LO_K \atop \vartriangleleft _H \in LO_H}
\Big\{ \la K ,\vartriangleleft _K\ra + \BM +\la H ,\vartriangleleft _H \ra,
\la H ,\vartriangleleft _H\ra ^*+ \BM ^* +\la K ,\vartriangleleft _K \ra^*\Big\} .
$$
In addition, since each element of $\L _{\Y _\X}^{F_{\bar c}} $ is a linear order without end points,
we have $K=H=\emptyset$ and, hence, $\BM=\BL _{\X _1}$ and $\L _{\Y _\X}^{F_{\bar c}}=\{ \BL _{\X _1} ,\BL _{\X _1} ^*\}$.
Thus each structure $\X _2 \in \CM ^{\bar c}_{\Y _\X}$ consists of the initial part, $\{c_0\}+\dots +\{c_{n-1}\}$,
labeled by the unary relations $U_j^{\X _2} =\{ c_j \}$, $j<n$,
and a final part,
which is either isomorphic to the linear order $\BL _{\X _1}$ or to its reverse, $\BL _{\X _1}^*$.
So we have $|(\CM ^{\bar c}_{\Y _\X})^{\cong }|\leq 2$.
\kdok

Finally we prove the second part of Theorem \ref{T030}.
By our analysis, the theory $\CT$ is $\o$-categorical iff Case A appears;
so, we have to prove that the $L_n$-structure $\X _0 =\la \o ,\prec , \{ a_0 \} ,\dots ,\{ a_{n-1 }\}\ra$ is $\o$-categorical iff
$\BL _{\X _0}:=\la \o \setminus F_{\bar a}, \prec \upharpoonright (\o \setminus F_{\bar a})\ra$ is an $\o$-categorical linear order.
Since $\Aut (\X _0)=\{ \id _{F_{\bar a}}\cup \,f : f\in \Aut (\BL _{\X _0}) \}$, for $n\in \N$ and
$\bar x ,\bar y \in (\o \setminus F_{\bar a})^n$ we have $\bar x \sim _{\BL _{\X _0}}\bar y \Leftrightarrow \bar x \sim _{\X _0}\bar y $,
which implies that $|(\o \setminus F_{\bar a})^n /\sim _{\BL _{\X _0}}|\leq|\o ^n /\sim _{\X _0}|$.
So, if $\X _0$ is $\o$-categorical, then $\BL _{\X _0}$ is $\o$-categorical (by the theorem of Engeler, Ryll-Nardzewski and Svenonius).

On the other hand, if $\BL _{\X _0}$ is $\o$-categorical, then the linear order $\la \o ,\prec \ra \cong n+ \BL _{\X _0}$
is $\o$-categorical (see Rosenstein's theorem, \cite{Rosen}, p.\ 299) and, since $\Aut (\X _0)=\Aut (\la \o ,\prec \ra )$,
$\X _0$ is $\o$-categorical too.
\kdok

\end{document}